\documentclass[12pt]{amsart}

\newcommand{\calC}{\mathcal{C}}
\newcommand{\calD}{\mathcal{D}}
\newcommand{\bbZ}{\mathbb{Z}}
\newcommand{\bbC}{\mathbb{C}}
\newcommand{\bbQ}{\mathbb{Q}}
\newcommand{\tensor}{\otimes}
\newcommand{\ext}[1]{\operatorname{\stackrel{#1}{\wedge}}}

\DeclareMathOperator{\Symp}{Sp}
\DeclareMathOperator{\HH}{H}
\theoremstyle{plain}
\newtheorem{lemma}{Lemma}
\newtheorem{thm}[lemma]{Theorem}

\begin{document}

\title{The Hodge Conjecture for general Prym varieties}

\author[I. Biswas]{Indranil Biswas}
\address{School of Maths, TIFR, Homi Bhabha Road, Mumbai 400 005, India.}
\email{indranil@math.tifr.res.in}

\author[K. H. Paranjape]{Kapil H.~Paranjape}
\address{IMSc, CIT Campus, Tharamani, Chennai 600 113, India.}
\email{kapil@imsc.ernet.in}

\maketitle

\section*{Introduction}

We work over $\bbC$, the field of complex numbers.

The Prym variety of a double cover $C\to D$ of a smooth connected
projective curve $D$ by a smooth connected curve $C$ is defined (see
\cite{Mumford}) as the identity component of the kernel of the norm
homomorphism $N:J(C)\to J(D)$ between the Jacobians of the curves.
This is an abelian variety polarised by the restriction of the
canonical polarisation on $J(C)$; we denote this variety by $P(C\to
D)$ or simply $P$ when there is no possibility of ambiguity.

A Hodge class on a variety $X$ is an integral singular cohomology
class on the complex manifold $X(\bbC)$ which is represented by a
closed differential form of type $(p,p)$. The Hodge conjecture (see
\cite{Hodge}) asserts that some multiple of such a class is the
cohomology class of an algebraic cycle on $X$.

Let $A$ be an abelian variety. The K\"unneth decomposition implies
that the rational singular cohomology of $A\times\dots\times A$ is a
direct sum of subquotients of tensor products of
$\HH^1(A(\bbC),\bbQ)$. Hence we have an action of a linear
automorphism of this vector space on these cohomology groups.  The
Mumford-Tate group $H(A)$ of $A$ can thus be defined (see \cite{DMOS})
as the group of all linear automorphisms of $\HH^1(A(\bbC),\bbQ)$
which stabilise all Hodge cycles on the varieties $A\times\dots\times
A$.

The aim of this note is to show that the Mumford-Tate group $H(P)$ of
a {\em general} Prym variety $P(C\to D)$ is isomorphic to the full
symplectic group $\Symp(2g)$; where the class in
$\ext{2}\HH^1(P(\bbC),\bbQ)=\HH^2(P(\bbC),\bbQ)$ which is stabilised
by this group is the first Chern class of the natural polarisation on
the Prym variety. Invariant theory (see \cite{Weyl} or \cite{Howe1}
and \cite{Howe2}) then implies that the only Hodge cycles on $P$ are
powers (under cup-product) of this polarisation class. In particular,
we obtain the Hodge conjecture for $P$ as a consequence of this
result.

As a particular case the N\'eron-Severi group of a general Prym
variety is $\bbZ$. This was proved earlier by Pirola (see \cite{Pirola}).
We do not give a new proof of that result and use it in an essential
way to prove our result.

The outline of the paper is as follows. In section~1 we set out some 
standard arguments about Mumford-Tate groups in families.
In section~2 we use an extension (due to Beauville \cite{Beauville})
of the definition of Prym varieties to the case where $C$ and $D$
are singular curves. The results on Mumford-Tate groups are applied
to this larger family of Prym varieties in section three. In addition
we use the semi-simplicity of the Mumford-Tate group (see \cite{DMOS})
and the result of Pirola (see \cite{Pirola}) to reduce the problem to an
elementary lemma on subgroups of the symplectic group.

\section{Mumford-Tate groups in families}

Let $f:X\to S$ be a family of smooth projective varieties parametrised
by a smooth connected variety $S$. For some positive integer $k$
let $V=R^if_*\bbQ_X$ denote the variation of pure Hodge structures of
weight $k$ on $S$. More generally we can consider any variation $V$
of Hodge structures  of weight $k$ on $S$.

Let $V^{a,b}= V^{\tensor a}\tensor V^{*\tensor b}$ be the associated
tensor variations of pure Hodge structures of weight $(a-b)k$. For
every $(a,b)$ such that $(a-b)k=2p$ is even, we have the nested sequence
of analytic subvarieties 
\[ H^{a,b}:=V^{a,b}_{\bbZ}\cap F^p V^{a,b} \subset V^{a,b}_{\bbZ}
                \subset V^{a,b}_{\bbC}
\]
of the complex vector bundle $V^{a,b}_{\bbC}$ over $S$ associated with
$V^{a,b}$.  The analytic variety $H^{a,b}$ parametrises pairs $(s,c)$,
where $s$ is a point of $S$ and $c$ an integral class of type $(p,p)$
in $V^{a,b}_s$; i.~e.\ $c$ is a {\em Hodge cycle}.

If $W$ is an irreducible component of $H^{a,b}$ such that the natural
map $W\to S$ is open at some point, then $W$ contains an open subset
of $V^{a,b}_{\bbZ}$; hence $W$ is a connected component of
$V^{a,b}_{\bbZ}$. Let $A^{a,b}$ be the the union of all such components.
The map $A^{a,b}\to S$ makes each component of the former a covering
space of $S$.

Now, if $W$ is an irreducible component of $H^{a,b}$ for which the
map $W\to S$ is {\em not} open at any point then its image in $S$ is
a set of measure zero by Sard's theorem. Let $B$ be the (countable)
union of these images as we vary over all the components of $H^{a,b}$
and as we vary $a$ and $b$.

If $s$ is any point of $S$ which is not in $B$, then by the above
reasoning, the only points of $H^{a,b}$ that lie over it are in
$A^{a,b}$.  Let $t$ be any other point of $S$ and $\gamma$ be a path in $S$
connecting $s$ and $t$. We can use $\gamma$ to identify $V^{a,b}_{\bbZ,s}$
with $V^{a,b}_{\bbZ,t}$; this then gives an identification of $A^{a,b}_s$
with $A^{a,b}_t$.  Hence, under this identification, the collection of
Hodge cycles in $V^{a,b}_{\bbZ,s}$ is contained in the collection of
Hodge cycles in $V^{a,b}_{\bbZ,t}$. Thus the Mumford-Tate group $G_t$
of $V_t$ is identified by $\gamma$ with a subgroup of the Mumford-Tate
group $G_s$ of $V_s$. In other words, we have
\begin{quote}
        The Mumford-Tate group at a general point contains (a
conjugate of) the Mumford-Tate group at a special point in a variation
of Hodge structures over a smooth connected variety. 
\end{quote}

\section{Degenerate covers}

A connected projective curve which has at worst ordinary double points
as its singularities is called a {\em semi-stable} curve.  The {\em
  dual graph} of such a curve has as its vertices the irreducible
components; each singular point gives an edge incident on the two
vertices corresponding to the components that contain it. We will be
interested in semi-stable curves whose dual graph is contractible and
hence a tree; such curves are called {\em tree-like}. In this case,
the first cohomology of the curve is a direct sum of the first
cohomology of its components with the induced (pure) Hodge structure.
In particular, the Jacobian of a tree-like semi-stable curve is the
product of the Jacobians of its components.

A finite morphism $C\to D$ of semi-stable curves is called a
semi-stable cover (or an admissible cover) if
\begin{enumerate}
\item This is a topological cover of constant degree
of $D$ outside a finite set of points which includes the singular
locus of $D$.
\item The inverse image of a singular point of $D$ consists of singular
points of $C$.
\item The order of ramification on the two branches
at a singular point of $C$ must be equal.
\end{enumerate}
This notion was first defined by Beauville \cite{Beauville} for the
case of degree two covers (which are the case of interest) and later
generalised (see \cite{Mumford-Harris}). In these papers, it is shown
that the deformations of such a semi-stable cover of tree-like curves
are unobstructed. In other words, there is a smooth (open) curve $S$,
a flat morphism $p:\calD\to S$ and a finite flat morphism
$f:\calC\to\calD$. There is a point $o$ of $S$ over which the $f$ restricts
to the given semi-stable cover $C\to D$. Moreover, the general
fibre is a double cover $C'\to D'$ of a smooth curve $D'$ by a smooth
curve $C'$.

We are interested in the case of degree two covers $C\to D$; here the
singular points of $C$ are either unramified on each branch or ramified
of order two on each branch.  Let us further assume that $C$ and $D$ are
tree-like.  For each component $D_i$ of $D$ there are two possibilities:
\begin{enumerate}
\item There is exactly one component $C_i$ of $C$ that lies over it.
The map $C_i\to D_i$ is a double cover in the usual sense.
\item There are two components $C'_i$ and $C''_i$ of $C$ that lie
over $D_i$ and the given map is an isomorphism between these
components and $D_i$.
\end{enumerate}
The Prym variety can be defined as before as the identity component of
the kernel of the natural norm homomorphism between the Jacobians
$J(C)\to J(D)$. It follows that the Prym variety is the product of the
Prym varieties of the covers $C_i\to D_i$ corresponding to the first
case and the Jacobians of the curves $D_i$ corresponding to the second
case. In particular, the product of these components gives an abelian
variety. Hence we have
\begin{quote}
  The family of Prym varieties can be extended to include the Prym
  varieties of degenerate tree-like covers. In particular, the Mumford-Tate
  group of a general Prym variety contains (a conjugate of) the
  Mumford-Tate group of the Prym variety of any degenerate tree-like cover.
\end{quote}

In the special case when $D$ has two exactly components (call them
$D_1$ and $D_2$), such a cover can be constructed in one of two ways:
\begin{enumerate}
\item[I] Let $C_1\to D_1$ be a double cover that is not branched at the
common point $p=D_1\cap D_2$. Then, $C$ is obtained by attaching to
$C_1$ two copies of $D_2$, one at each point lying over $p$.
\item[II] Let $C_1\to D_1$ and $C_2\to D_2$ be double covers that are
both branched at the common point $p$. We obtain $C$ by attaching
the curves $C_1$ and $C_2$ along their respective ramification points
lying over $p$.
\end{enumerate}
The specific covers that we are interested in are the following.
\begin{enumerate}
\item A covering of type (II) which is the degeneration of a double
  cover $C\to D$ where $D$ is rational and $C$ is of genus $g$.  The
  curves $D_i$ are smooth rational curves. The curve $C_1$ is an
  elliptic curve and the curve $C_2$ is any hyperelliptic curve of
  genus $g-1$.
\item A covering of type (I) which is the degeneration of a double
  cover $C\to D$ where $D$ has genus at least 2 and the cover is
  \'etale. The curve $D_1$ is any elliptic curve, $C_1\to D_1$ is an
  \'etale double cover and $D_2$ is any curve of genus one less than
  that of $D$.
\item A covering of type (I) which is the degeneration of a double
  cover $C\to D$ where $D$ has genus at least 1 and the cover is
  ramified at some point. The curve $D_1$ is any curve of genus one
  less than that of $D$ and $C_1\to D_1$ is a double cover ramified at
  the same number of points as the cover $C\to D$; $D_2$ is any
  elliptic curve.
\end{enumerate}
As a result we have
\begin{lemma}\label{main} We have containments of Mumford-Tate groups
  as enumerated below.
\begin{enumerate}
\item The Mumford-Tate group of a general hyperelliptic curve of genus
  $g$ contains a conjugate of the product of the Mumford-Tate group of
  any elliptic curve with the Mumford-Tate group of any hyperelliptic
  curve of genus $g-1$.
\item The Mumford-Tate group of the Prym variety of a general \'etale
  cover of a curve of genus $g\geq 2$ contains a conjugate of the
  Mumford-Tate group of any curve of genus $g-1$.
\item The Mumford-Tate group of the Prym variety of a general cover of
  a curve of genus $\geq 1$ ramified at $r\geq 1$ points contains a
  conjugate of the product of the Mumford-Tate group of any elliptic
  curve with the Mumford-Tate group of the Prym variety of any cover
  of a curve of genus $g-1$ which is ramified at $r$ points.
\end{enumerate}
\end{lemma}
\begin{proof}
  The first cohomology group of the product of two abelian varieties is
  the direct sum of the first cohomology groups of the individual
  abelian varieties. Moreover, the Hodge cycles on the individual
  varieties pull-back to give Hodge cycles on the product. Thus it
  follows that the Mumford-Tate group of the product contains the
  product of the Mumford-Tate groups. The result now follows from the
  above constructions.
\end{proof}

\section{The Main result}

To prove the main result we need the following three lemmas.
\begin{lemma}[Pirola] The N\'eron-Severi group of a general Prym
  variety is the free group on 1 generator. 
\end{lemma}
This lemma is proved in \cite{Pirola}. We note that this case includes
the case of a general hyperelliptic curve.
\begin{lemma} Let $G$ be a connected semi-simple subgroup of the
  symplectic group $Sp(2n)$ which contains (a conjugate of) the
  product $Sp(2a)\times Sp(2n-2a)$, then $G$ is either this product or
  it is the the full symplectic group
\end{lemma}
\begin{proof}
  Let $V$ be the standard representation of $Sp(2n)$. Let
  $\oplus_{i\in I} W_i$ be its decomposition into isotypical
  components as a representation of $G$. Let $V=V_1\oplus V_2$ be the
  decomposition of $V$ as a representation of $Sp(2a)\times
  Sp(2n-2a)$. Then each $W_i$ is either $V_1$ or $V_2$ or $V=V_1\oplus
  V_2$. The result follows by dimension counting.
  
  The lemma also follows from the fact that the quotient
  \[ \frac{sp(2n)}{sp(2a)\times sp(2n-2a)} 
  \]
  of the Lie algebras is an irreducible module over $Sp(2a)\times
  Sp(2n-2a)$.
\end{proof}
\begin{lemma}
  Let $A$ be an abelian variety of dimension $n$ whose Mumford-Tate
  group is the product $Sp(2a)\times Sp(2n-2a)$ in $Sp(2n)$. Then the
  N\'eron-Severi group of $A$ is of rank at least 2.
\end{lemma}
\begin{proof}
  The first cohomology group of $A$ decomposes as a direct sum of two
  (polarised) sub-Hodge structures. It follows that $A$ is the product
  of two abelian subvarieties. Hence we have the result.
\end{proof}
\begin{thm}\label{Main Theorem}
  The Mumford-Tate group of a general Prym variety is the full
  symplectic group.
\end{thm}
\begin{proof}
  We begin with the case where the base curve has genus zero. In this
  case the Prym varieties are the Jacobians of the corresponding
  hyperelliptic double cover. The result is classical for elliptic
  curves which can be considered as the Prym varieties associated with
  double covers of smooth rational curves branched at 4 points. By
  induction, let us assume that the result is known for hyperelliptic
  Jacobians of genus less than $g\geq 2$. The lemma~\ref{main} then
  shows that the Mumford-Tate group of a general hyperelliptic curve
  of genus $g$ contains $Sp(2)\times Sp(2g-2)$. By the above results
  we see that thus Mumford-Tate group must be either $Sp(2g)$ or
  $Sp(2)\times Sp(2g-2)$. In the latter case, the N\'eron-Severi group
  of the curve would have rank at least two but by Pirola's result we
  know that this is not true for the general hyperelliptic curve.
  Hence we see that the Mumford-Tate group of a general hyperelliptic
  curve must be $Sp(2g)$ where $g$ is the genus of the curve.
  
  Now let us consider the case where the cover is unramified. Then we
  may assume that the base curve of the double cover has genus $g$ at
  least 2 (else the Prym variety is just a point). In this case the
  Prym variety has dimension $n=g-1$. By the lemma~\ref{main} we know
  that the Mumford-Tate group contains the Mumford-Tate group of any
  curve of genus $g-1$. In particular, it contains the Mumford-Tate
  group of any hyperelliptic curve of genus $g-1$ and hence by the
  previous paragraph it contains (and is thus equal to) $Sp(2n)$.
  
  Now assume that the base curve of the double cover has genus $g$ at
  least 1 and the cover is ramified. We argue by induction on the
  genus of the base curve. We can begin the induction since we already
  know the result for the hyperelliptic curves. Let us assume that the
  result is known for base curves of genus less than $g$. By the
  lemma~\ref{main} we know that the Mumford-Tate group contains the
  product of the Mumford-Tate group of an elliptic curve with the
  Mumford-Tate group of the Prym variety of a the double of a curve of
  genus $g-1$; in other words it contains $Sp(2)\times Sp(2n-2)$ by
  induction. Now as argued above, the three lemmas above imply that
  the Mumford-Tate group must be $Sp(2n)$.
\end{proof}


\providecommand{\bysame}{\leavevmode\hbox to3em{\hrulefill}\thinspace}

\end{document}